   \documentclass[11pt]{amsart}
   \usepackage{amsmath,amssymb}
   
   \renewcommand{\bf}{\bfseries}

   \newcommand{\C}{\mathbb{C}}
   \newcommand{\E}{\mathbb{E}}
   \newcommand{\N}{\mathbb{N}}
   \newcommand{\Z}{\mathbb{Z}}
   
   \newtheorem{theorem}{Theorem}
   
   \newtheorem{corollary}{Corollary}
   \newtheorem{lemma}{Lemma}
   \newcommand{\dis}{\displaystyle}

\begin{document}
   \today

   \title{Pointwise convergence along cubes for measure preserving systems}
   \author{by I. Assani}
   \thanks{Department of Mathematics,
       UNC Chapel Hill, NC 27599\\AMS subject classification  37A30 28D05 \\email: assani@email.unc.edu}

   \begin{abstract}
   Let $(X, \mathcal{B}, \mu)$ be a probability measure space and
   $T_1$, $T_2$ , $T_3$ three not necessarily commuting measure
   preserving transformations on $(X, \mathcal{B}, \mu)$.
   We prove that for all bounded functions $f_1$, $f_2$, $f_3$ the
   averages
   $$\frac{1}{N^2}\sum_{n, m =1}^N
   f_1(T_1^nx)f_2(T_2^mx)f_3(T_3^{n+m}x)$$ converges a.e..
   Generalizations to averages of $2^k -1$ functions are also
   given for not necessarily commuting weakly mixing systems.
  \end{abstract}
  \maketitle
  \section{Introduction}
  In [A1] and [A2] we proved that if $T$ is a measure preserving
  transformation on $(X, \mathcal{B},\mu)$ then the averages
  of three functions
  $$\frac{1}{N^2}\sum_{n, m=1}^N f_1(T^nx)f_2(T^mx)f_3(T^{n+m}x)$$ or more generally $2^k -1$ functions converge
  a.e.

 We want to show that the method we used in these papers can yield  more general pointwise
 results. More precisely we want to show that one can have pointwise convergence when $T$ is replaced by measure preserving
 transformations $T_i$, $1\leq i\leq 3$ that do not necessarily commute. As shown in [Be]
 Khintchin 's recurrence theorem [Kh] can be extended by the
 convergence of such averages. One can observe that if $T_1$ and $T_2$ do
not necessarily commute then the averages
$$\frac{1}{N}\sum_{n=1}f(T_1^nx)g(T_2^nx)$$ may diverge ([Ber]).
Also an example given in [L] shows that the averages
$$\frac{1}{N^2}\sum_{n,m=1}^N\mu(A\cap T_1^{-n}A\cap T_2^{-m}A\cap
T_1^{-n}T_2^{-m}A)$$ may also diverge if $T_1$ and $T_2$ do not
necessarily commute.

 \begin{theorem}
  Let $(X, \mathcal{B}, \mu)$ be a probability measure space and
  $T_1$, $T_2$ , $T_3$ three not necessarily commuting measure
  preserving transformations on $(X, \mathcal{B}, \mu)$.
  Then for all bounded functions $f_i$, $1\leq i \leq 3$ the
  averages
   $$\frac{1}{N^2}\sum_{n, m =1}^N
   f_1(T_1^nx)f_2(T_2^mx)f_3(T_3^{n+m}x)$$ converge a.e.
\end{theorem}

 At the present time we do not know if the pointwise convergence
 holds for averages along the cubes of $2^k -1 $ functions for
 $k>2$ for not necessarily commuting measure preserving
 transformations.
 However if the transformations $T_i$, $1\leq i\leq k$ are weakly mixing then we can establish the
 pointwise convergence of the averages for all positive integer
 $k$ and identify the limit.
  \begin{theorem}
  Let $(X, \mathcal{B}, \mu)$ be a probability measure space and
 $T_i$ weakly mixing transformations (not necessarily commuting) on this measure space.
 Then the averages along the cubes applied to the bounded
 functions $f_i$, $1\leq i\leq 2^k -1$ converge a.e. to $\prod_{i=1}^{2^k-1} \int
 f_id\mu$.
 \end{theorem}
  The norm convergence follows by
 integration as the functions are in $L^{\infty}$.
 We can derive the following corollaries. The first one extends
 Khintchine's recurrence theorem. The case $T_1= T_2 = T$ was
 treated in [Be].
 \begin{corollary}
Let $(X, \mathcal{F}, \mu)$ be a probability measure space and
$T_1$, $T_2$ two measure preserving transformations on this
measure space. We denote by $\mathcal{I}_1$ and $\mathcal{I}_2$
the $\sigma$ algebras of the invariant sets for $T_1$ and $T_2$.
Consider $A$ a set of positive measure. Then
$$\lim_N \frac{1}{N^2}\sum_{n, m=1}^N \mu(A\cap T_1^{-n}A\cap
T_2^{-n-m}A)= \int_A\E(\mathbf{1}_A,
\mathcal{I}_1)(x).\E(\mathbf{1}_A, \mathcal{I}_2)(x)d\mu.$$ In
particular if $\mathcal{I}_1\subset \mathcal{I}_2$ (or
$\mathcal{I}_2\subset \mathcal{I}_1$) then
$$\lim_N \frac{1}{N^2}\sum_{n, m=1}^N \mu(A\cap T_1^{-n}A\cap
T_2^{-n-m}A)\geq \mu(A)^3.$$
\end{corollary}
 The assumption $\mathcal{I}_1\subset \mathcal{I}_2$ is satisfied
 if $T_1$ is ergodic as  the invariant functions for $T_1$ are then the
 constant functions.
 \smallskip

 We recall that a set of integers is said to be syndetic (also called relatively dense) if it has bounded
 gaps. A corollary of theorem 2 is the following.
 \begin{corollary}
   Let $(X, \mathcal{B}, \mu)$ be a probability measure space and
  $T_i$ weakly mixing transformations (not necessarily commuting) on this measure
  space and $0\leq \lambda<1$. For all measurable set A of positive measure, for all $k\geq 1$,
  for $\mu$ a.e. $x$ the set
  $$ \{(n_1,n_2,...,n_k)\in\Z^k: \mathbf1_A(x).\mathbf1_A(T^{n_1}_{1}x).\mathbf1_A(
  T^{n_1+n_2}_{2}x)\cdots\mathbf1_A(T^{n_1+n_2 +...+n_k}_{k}x)>\lambda \mu(A)^{2^k}\}$$
  is syndetic.
\end{corollary}

\section{Proof of theorem 1}
 The following lemma will be useful for the theorems we want to
 prove.
 \begin{lemma}
 Let $a_n$, $b_n$ and $c_n$, $n\in \N$ be three sequences of
  scalars that we assume for simplicity bounded by one. Then for each $N$ positive integer we have
 $$
\big|\frac{1}{N^2}\sum_{m,n=1}^N a_n. b_{m}. c_{n+m}\big|^2 \leq
4Min\bigg[\sup_t\big|\frac{1}{2N}\sum_{m'=1}^{2N}c_{m'}e^{2\pi
 im't}\big|^2, \sup_t\big|\frac{1}{N}\sum_{n'=1}^Na_{n'}e^{2\pi
in't}\big|^2\bigg]
$$
\end{lemma}
\begin{proof}
We denote by $M_N(a, b, c)$ the quantity
$\dis\frac{1}{N^2}\sum_{n, m=1}^N a_n.b_m.c_{n+m}$. The steps are
similar to those given in the proof of theorem 4 in [A1]. We have
\[
\begin{aligned}
&|M_N(a, b, c)|^2 \\
&\leq
\|a\|_{\infty}^2\bigg(\frac{1}{N}\sum_{n=1}^{N}\big|\frac{1}{N}\sum_{m=1}^{N}b_{m}c_{n+m}\big|^2\bigg)
\\
&\leq \|a\|_{\infty}^2 \frac{1}{N}\sum_{n=1}^{N}\bigg|\int
\big(\sum_{m=1}^{N}b_m e^{-2\pi
imt}\big)\big(\frac{1}{N}\sum_{m'=1}^{2N}c_{m'}e^{2\pi im't}\big).
e^{2\pi int} dt \bigg|^2 \\
&\leq \|a\|_{\infty}^2\frac{1}{N}\int \bigg|
\sum_{m=1}^{N}b_{m}e^{-2\pi
imt}\bigg|^2\bigg|\frac{1}{N}\sum_{m'=1}^{2N}c_{m'}e^{2\pi
im't}\bigg|^2dt \\
& \leq
4\frac{\|a\|_{\infty}^2}{N}\sup_t\bigg|\frac{1}{2N}\sum_{m'=1}^{2N}c_{m'}e^{2\pi
im't}\bigg|^2\int \big|\sum_{m=1}^{N}b_{m}e^{-2\pi
imt}|^2dt \\
& \leq
4\|a\|_{\infty}^2\sup_t\bigg|\frac{1}{2N}\sum_{m'=1}^{2N}c_{m'}e^{2\pi
im't}\bigg|^2\frac{1}{N} N\|b\|_{\infty}^2\\
& \leq
4\|a\|_{\infty}^2\|b\|_{\infty}^2\sup_t\big|\frac{1}{2N}\sum_{m'=1}^{2N}c_{m'}e^{2\pi
 im't}\big|^2
\end{aligned}
\]
 This provides a first bound for $|M_N(a, b, c)|^2$. To obtain the
 second bound we can start instead in the following manner.
 \[\begin{aligned}
 &|M_N(a, b, c)|^2 \\
 &\leq
 \|b\|_{\infty}^2 \frac{1}{N}\sum_{m=1}^N \bigg|\int
 \big(\frac{1}{N}\sum_{n=1}^N b_ne^{-2\pi
 int}\big)\big(\sum_{n'=1}^{2N}c_{n'}e^{2\pi in't}\big)e^{2\pi
 mt}dt\big)\bigg|^2
 \end{aligned}
\]
 From these last steps by using a similar path we obtain the second
 bound.
\end{proof}

  The Wiener-Wintner pointwise ergodic theorem asserts that if $T$ is a measure preserving transformation
on the probability measure space $(X, \mathcal{B}, \mu)$ and $f$ a
$L^{\infty}$ function then we can find a set of full measure $X_f$
such that for x in this set the averages
 \begin{equation}
 \frac{1}{N}\sum_{n=1}^N f(T^nx)e^{2\pi int}
\end{equation}
 converge for all real number t.
 One can see [A3], for instance, for various proofs of this result. The following lemma extends this result.

 \begin{lemma}
 Let $T_2$ and $T_3$ be two measure preserving transformations on
 $(X, \mathcal{B}, \mu)$. For each pair of functions $f_2$, $f_3$ in $L^{\infty}$ there exists a set of full measure
 $X_{f_2, f_3}$ such that if $x$ is in this set then the averages
 $$\frac{1}{N^2}\sum_{m,n=1}^N f_2(T_2^mx)f_3(T_3^{m+n}x)e^{2\pi
 int}$$ converge for all t.
 \end{lemma}
 \begin{proof}

 Without loss of generality we can assume that the functions $f_2$
 and $f_3$ are bounded by 1.
 \vskip1ex
 We consider an ergodic decomposition $\mu_{c, 3}$ for $T_3$ on $(X,
 \mathcal{B}, \mu)$. This means that on $(X, \mathcal{B}, \mu_{c,
 3})$ the transformation $T_3$ is measure preserving and ergodic.
  Furthermore $\mu_{c, 3}$ is a disintegration of $\mu$,
 i.e. for each integrable function $f\in L^1(\mu)$ we have
  $\int f(x)d\mu(x) = \int f(y)d{\mu}_{c}(y)dP(c)$ where P is a
  probability measure.
  \vskip1ex
   Using this ergodic decomposition we can conclude that for P
   a.e. $c$, for each positive integer $m$ the functions $f_2\circ T_2^m$ are all
   in $L^{\infty}(\mu_{c, 3})$ and bounded by one.
  The functions $f_3\circ T_3^m$ are also for P a.e. c in
  $L^{\infty}(\mu_{c,3})$.  So we consider the set $\overline{C_{3,1}}$
  of full measure where all these functions are bounded by one
  for $\mu_{c, 3}$ a.e. $y$. We restrict this set further by considering the disintegration of the set of $x$ where
  the averages
  \begin{equation}
  \frac{1}{N}\sum_{m=1}^N f_2(T_2^mx)e^{2\pi im\epsilon}
  \end{equation}
   converge for all $\epsilon$. This means that for P a.e. c there exists a set
   a set of $\mu_{c,3}$ full measure such that the averages in (2) converge
    for all $\epsilon$ real. Let us denote by $\overline{C_{3,2}}$ this set of full measure of c. Now we pick $c$ in the set
    $\dis \overline{C_3}= \overline{C_{3,1}}\cap \overline{C_{3,2}}$ and restrict ourselves to $(X, \mathcal{B}, \mu_{c, 3})$.
  We denote by $\mathcal{K}_{c, 3}$ the Kronecker factor of $T_3$. It consists of the closed linear span of the eigenfunctions
   of $T_3$ in $L^2(\mu_{c, 3})$ with an orthonormal bassis
   $e_{c,3}^k$ of eigenfunctions with modulus 1.
  We decompose the function $f_3$ into the sum $P_{\mathcal{K}_{c,
  3}}(f_3) + f - P_{\mathcal{K}_{c,3}}(f_3)$. The function $g_{c,
  3}= f - P_{\mathcal{K}_{c,3}}(f_3)$ being in the orthogonal
  complement of $\mathcal{K}_{c, 3}$ we have by the uniform Wiener
  Wintner ergodic theorem (see [A3]) for instance ) for $\mu_{c,
  3}$ a.e $y$
  \begin{equation}
  \lim_N\sup_t \bigg|\frac{1}{N}\sum_{m=1}^N g_{c,3}(T_3^{m}y)e^{2\pi
   imt}\bigg|=0.
  \end{equation}
 Applying lemma 1 pointwise with $a_n = e^{2\pi int}$, $b_m =
 f_2(T_2^my)$ and $c_{n+m}= g_{c,3}(T^{n+m}y)$ and using (3)  we obtain

 $$\lim_N\sup_t\bigg|\frac{1}{N^2}\sum_{m,n=1}^N f_2(T_2^my)g_{c,3}(T_3^{m+n}y)e^{2\pi
 int}\bigg|=0.$$
It remains to prove the convergence of
$$\frac{1}{N^2}\sum_{m,n=1}^N f_2(T_2^my)P_{\mathcal{K}_{c,
  3}}(f_3)(T_3^{m+n}y)e^{2\pi
 int}$$ for all t.
 The function $P_{\mathcal{K}_{c,
  3}}(f_3)$ can be written in terms of the orthonormal basis
  $e_{c,3}^k$ as $\dis \sum_{k=1}^{\infty}\big(\int
  f_3\overline{e_{c,3}^k} d\mu_{c, 3}(y)\big).e_{c,3}^k.$
  For each eigenfunction $e_{c,3}^k$ with eigenvalue $\lambda_{c,k}$
  we have
  $$\frac{1}{N^2}\sum_{m,n=1}^N f_2(T_2^my)e_{c,3}^k(T_3^{m+n}y)e^{2\pi
 int}= e_{c, 3}(y)\frac{1}{N^2}\sum_{m, n=1}^N f_2(T_2^my)e^{2\pi
 i(m+n)\lambda_{c,k}}e^{2\pi int}.$$
The last term is equal to $\dis e_{c,3}(y)\frac{1}{N}\sum_{n=1}^N
e^{2\pi in(t + \lambda_{c, k})}\frac{1}{N}\sum_{m=1}^N
f_2(T_2^my)e^{2\pi im\lambda_{c,k}}.$ \vskip1ex

The sequence $\dis e_{c,3}(y)\frac{1}{N}\sum_{n=1}^N e^{2\pi in(t
+ \lambda_{c, k})}$ converges for all t by the convergence of
$\dis \frac{1}{N}\sum_{n=1}^n e^{2\pi in \theta}$ for each
$\theta$ real. The Wiener Wintner ergodic theorem and the
disintegration mentioned above guarantee the convergence of $\dis
\frac{1}{N}\sum_{m=1}^N f_2(T_2^my)e^{2\pi im\lambda_{c,k}}$ for
$\mu_{c,3}$ a.e. $y$. By linearity we can reach the same
conclusion for the finite sum $\dis \sum_{k=1}^K\big(\int
  f_3.\overline{e_{c,3}^k} d\mu_{c, 3}(y)\big).e_{c,3}^k.$
  The same conclusion for
  $\dis P_{\mathcal{K}_{c,
  3}}(f_3) = \sum_{k=1}^{\infty}\big(\int
  f_3\overline{e_{c,3}^k} d\mu_{c, 3}(y)\big).e_{c,3}^k $

  follows by approximation and the use of the maximal inequality
  in $L^2(\mu_{c,3}).$

  Thus we have found a set of c of full P measure such that for $\mu_{c, 3}$ a.e, $y$ the averages
$$\frac{1}{N^2}\sum_{m,n=1}^N f_2(T_2^my)f_3(T_3^{m+n}y)e^{2\pi
 int}$$ converge for all t.
 By integrating with respect to c we obtain a set of x of full
 measure for $\mu$ where
 $$\frac{1}{N^2}\sum_{m,n=1}^N f_2(T_2^mx)f_3(T_3^{m+n}x)e^{2\pi
 int}$$ converge for all t. This concludes the proof of the lemma.
\end{proof}
\vskip1ex
 \noindent{\bf End of the proof of theorem 1}
 With the previous lemmas we can finish the proof of theorem 1.
 We take an ergodic decomposition of $T_1$ with respect to $\mu$
 We denote the disintegrated measures by $\mu_{c,1}$.
 By using the previous lemma for $f_2$ and $f_3$ fixed functions in $L^{\infty}(\mu)$ we can
 find a set of full measure $\overline{D}$ such that if c is this
 set then we have the following properties;
\begin{enumerate}
\item the functions $f_1\circ T_1^n(y)$, $f_2\circ T_2^m(y)$ and
$f_3\circ T_3^{n+m}(y)$ are $\mu_{c,1}$ a.e. $y$ bounded by one

\item for $\mu_{c, 1}$ a.e. $y$ the sequence $\dis
\frac{1}{N^2}\sum_{m,n=1}^N f_2(T_2^my)f_3(T_3^{m+n}y)e^{2\pi
 int}$ converges for all real number $t$.
 \end{enumerate}
 \vskip1ex

 We fix $c$ in $\overline{D}$ and denote by $\mathcal{K}_{c, 1}$ the Kronecker factor of $T_1$.
 We decompose the function $f_1$ into the sum $\dis P_{\mathcal{K}_{c,
  1}}(f_1) + f - P_{\mathcal{K}_{c,1}}(f_1).$
   The function $P_{\mathcal{K}_{c,1}}(f_1)$ can be written as
   $\dis \sum_{k=1}^{\infty}\big(\int
  f_1.\overline{e_{c,1}^k} d\mu_{c, 1}(y)\big).e_{c,1}^k$ where the
  functions $e_{c,1}^k$ are eigenfunctions for $T_1$ of modulus
  one with eigenvalues $\alpha_{c,k}.$
  We can use (2) above to prove the convergence of the averages
  $$\frac{1}{N^2}\sum_{m,n=1}^Ne_{c,1}^k(T_1^ny)f_2(T_2^my)f_3(T_3^{m+n}y).$$
  By linearity and approximation we can prove the convergence for $\mu_{c,1}$ a.e. $y$ of the
  averages
  $$\frac{1}{N^2}\sum_{n, m=1}^N P_{\mathcal{K}_{c,
  1}}(f_1)(T_1^ny)f_2(T_2^my)f_3(T_3^{m+n}y).$$
  The convergence of the averages
$$\frac{1}{N^2}\sum_{n, m=1}^N [f_1-P_{\mathcal{K}_{c,
  1}}(f_1)](T_1^ny)f_2(T_2^my)f_3(T_3^{m+n}y)$$
  is obtained by applying pointwise the second bound listed in
  lemma 1. We pick $a_n= [f_1-P_{\mathcal{K}_{c,
  1}}(f_1)](T_1^ny)$, $b_m= f_2(T_2^my)$ and $c_{n+m}=
  f_3(T_3^{n+m}y)$. The result follows by the uniform Wiener
  Wintner theorem applied to the function $[f_1-P_{\mathcal{K}_{c,
  1}}(f_1)]$ and the ergodic dynamical system $(X, \mathcal{B},
  \mu_{c,1}, T_1)$.
  We can finish the proof by integrating with respect to $P.$
\section{Proof of Theorem 2}
 The proof can be made by induction on $k$.
 \vskip1ex
 \noindent{\bf The case k=2}
\vskip1ex

 We have in this case the following lemma.
 \begin{lemma}
 Let $(X, \mathcal{B}, \mu)$ be a probability measure space and
 $T_1$, $T_2$ and $T_3$ be three weakly mixing measure preserving transformations on
 this space. Then for all $L^{\infty}$ functions, $f_1$, $f_2$
 and $f_3$ the averages
 $$\frac{1}{N^2}\sum_{m, n =1}^N
 f_1(T_1^{n}x)f_2(T_2^{m}x)f_3(T_3^{n+m}x)$$ converge a.e. to
 $\prod_{i=1}^3 \int f_i d\mu.$
 \end{lemma}
 \begin{proof}
 The lemma follows from the proof of theorem 1. When the transformations are weakly mixing the Kronecker factors
 are all reduced to the constant functions identified with $\C$.
 Thus the pointwise limit will be zero for $\mu$ a.e. $x$  if one of the functions $f_i$, $1\leq i\leq 3$ has zero
 integral. The result follows without difficulty from this observation.
 \end{proof}
 \vskip1ex
\noindent{\bf The case $k>2$}
\vskip1ex

 The induction method will be sufficiently described by considering the case
 $k=3$. Moving to higher values of $k$ can be done in the same way
 as in [A1]. We only sketch the proof as we can follow a similar
 path.

 So we consider seven weakly mixing transformations on $(X,
 \mathcal{B}, \mu)$, $T_i$, $1\leq i\leq 7$ and seven bounded
 functions $f_i$, $1\leq i\leq 7$. For simplicity we denote
  $f(T^mx)$ by $T^mf(x).$
The averages in this case are
\[
\begin{aligned}
&M_N(f_1, f_2,\cdots ,f_7)(x)\\
&= \frac{1}{N^3}\sum_{n, m, p=1}^N
T_1^nf_1(x)T_2^nf_2T_3^pf_3(x)T_4^{n+m}f_4(x)T_5^{n+p}f_5(x)T_6^{p+m}f_6(x)T_7^{n+m+p}f_7(x)
\end{aligned}
\]
We have the following lemma.
\begin{lemma}
If $f_1$ or $f_2$ is in $\C^{\perp}$ then for a.e. $x$
\begin{equation}
\lim_N\frac{1}{N}\sum_{n=1}^{N}\sup_t\bigg|\frac{1}{N}\sum_{m=1}^{N}T_1^mf_1(x)T_2^{n+m}f_2(x)e^{2\pi
imt}\bigg|^2 = 0
\end{equation}
\end{lemma}
\begin{proof}
This can be obtained by following the same steps as those used in
[A1]. The assumption made that $f_1$ or $f_2$ are in $\C^{\perp}$
is reflected in the fact that $\dis \lim_H\frac{1}{H}\sum_{h=1}^H
\big|\int T_1^nf_1T_1^{n+h}f_1 d\mu\big|= 0.$ (one can assume that
the functions are real). We skip the proof of this lemma.
\end{proof}
\noindent{\bf End of the proof of theorem 2}
\vskip1ex
\[
\begin{aligned}
&|M_N(f_1,f_2,...,f_7)|^2 \\
&=\bigg|
\frac{1}{N^3}\sum_{p=1}^{N}T_1^pf_1(x)\sum_{n=1}^{N}T_2^nf_2(x)T_3^{p+n}f_3(x)\big(\sum_{m=1}^{N}T_4^mf_4(x)T_5^{n+m}f_5(x)T_6^{p+m}f_6(x)
T_7^{n+m+p}f_7(x)\big)\bigg|^2 \\
&\leq
\frac{1}{N^2}\sum_{p=1}^{N}\sum_{n=1}^{N}\|f_1\|_{\infty}^2\|f_2\|_{\infty}^2\|f_3\|_{\infty}^2\bigg|\frac{1}{N}
\sum_{m=1}^{N}T_4^mf_4(x)T_5^{n+m}f_5(x)T_6^{p+m}f_6(x)T_7^{p+n+m}f_7(x)\bigg|^2
\\
&= \frac{1}{N^2}\prod_{i=1}^{3}\|f_i\|_{\infty}^2. \\
&\sum_{n=1}^{N}\sum_{p=1}^{N}\bigg|\int \big( \sum_{m=1}^{(N)}
T_4^mf_4(x)T_5^{n+m}f_5(x)e^{-2\pi
imt}\big)\big(\frac{1}{N}\sum_{m'=1}^{2N}T_6^{m'}f_6(x)T_7^{n+m'}f_7(x)e^{2\pi
im't}\big). e^{2\pi ipt}dt\bigg|^2 \\
&\leq
\frac{1}{N^2}\prod_{i=1}^3\|f_i\|_{\infty}^2\sum_{n=1}^{N}\int
\bigg| \sum_{m=1}^{N} T_4^mf_4(x)T_5^{n+m}f_5(x)e^{-2\pi
imt}\big)\big(\frac{1}{N}\sum_{m'=1}^{2N}T_6^{m'}f_6(x)T_7^{n+m'}f_7(x)e^{2\pi
im't}\big)\bigg|^2 dt \\
&\leq
\frac{C}{N^2}\prod_{i=1}^3\|f_i\|_{\infty}^2\sum_{n=1}^{N}\sup_t
\bigg|\frac{1}{N}\sum_{m'=1}^{N}T_6^{m'}f_6(x)T_7^{n+m'}f_7(x)e^{2\pi
im't}\bigg|^2 N\prod_{j=4}^{5}\|f_j\|_{\infty}^2 \\
&= C\prod_{i=1}^{5}\|f_i\|_{\infty}^2
\frac{1}{N}\sum_{n=1}^{N}\sup_t
\bigg|\frac{1}{N}\sum_{m'=1}^{N}T_6^{m'}f_6(x)T_7^{n+m'}f_7(x)e^{2\pi
im't}\bigg|^2
\end{aligned}
\]
With the help of lemma 4 one can conclude that if $f_6$ or $f_7$
belong to $\C^{\perp}$ then the averages of these seven functions
converge to zero. By using the symmetry of the sum of the averages
with respect to $n$, $m$ and $p$ one can see that the averages
will converge to zero if one of the functions $f_i\in \C^{\perp},$
$1\leq i\leq 7$.
\section{Proof of the Corollaries}
\subsection{Corollary 1}
 The averages
 $$\frac{1}{N^2}\sum_{n, m=1}^N \mu(A\cap T_1^{-n}A\cap
T_2^{-n-m}A)$$ are the integrals of the functions
$$\frac{1}{N^2}\sum_{n, m = 1}^N \mathbf{1}_A(x)\mathbf{1}_A
(T_1^n x)\mathbf{1}_A(T_2^{n+m}x)$$ with respect to the measure
$\mu$. As a particular case of theorem 1 we have the pointwise
convergence of these averages. Thus $$\lim_N\frac{1}{N^2}\sum_{n,
m=1}^N \mu(A\cap T_1^{-n}A\cap T_2^{-n-m}A)$$ exists after
integration. So we just have to prove that
$$\lim_N\frac{1}{N^2}\sum_{n, m=1}^N
\mathbf{1}_A(T_1^nx)\mathbf{1}_A(T_2^{n+m}x)= \E(\mathbf{1}_A,
\mathcal{I}_1)(x).\E(\mathbf{1}_A, \mathcal{I}_2)(x)$$ in $L^2$
norm to conclude. For each $N$ we have
\[\begin{aligned}
&\frac{1}{N^2}\sum_{n, m=1}^N
\mathbf{1}_A(T_1^nx)\mathbf{1}_A(T_2^{n+m}x)\\
&= \frac{1}{N^2}\sum_{n,m=1}^N
\mathbf{1}_A(T_1^nx)\E(\mathbf{1}_A, \mathcal{I}_2)(x) +
\frac{1}{N^2}\sum_{n, m=1}^N
\mathbf{1}_A(T_1^nx)[\mathbf{1}_A(T_2^{n+m}x)- \E(\mathbf{1}_A,
\mathcal{I}_2)(x)]
\end{aligned}
\]
The first term of the last equation converges by Birkhoff's
pointwise ergodic theorem to $\E(\mathbf{1}_A,
\mathcal{I}_1)(x).\E(\mathbf{1}_A, \mathcal{I}_2)(x)$. Noticing
that the function $\E(\mathbf{1}_A, \mathcal{I}_2)(x)$ is $T_2$
invariant we can bound the $L^2$ norm of the second term by
 $$\|\frac{1}{N}\sum_{n=1}^N\big|\frac{1}{N}\sum_{m=1}^N[\mathbf{1}_A\circ
 T_2^m - \E(\mathbf{1}_A, \mathcal{I}_2)]\circ T_2^n\|_2.$$
 This term is less than
 $$\frac{1}{N}\sum_{n=1}^N \|\sum_{n=1}^N\big|\frac{1}{N}\sum_{m=1}^N[\mathbf{1}_A\circ
 T_2^m - \E(\mathbf{1}_A, \mathcal{I}_2)]\|_2$$
 which is equal to
 $$\|\frac{1}{N}\sum_{m=1}^N[\mathbf{1}_A\circ
 T_2^m - \E(\mathbf{1}_A, \mathcal{I}_2)]\|_2$$
 This last term tends to zero by the mean ergodic theorem applied
 to $T_2$. This proves that
 $\dis \lim_N\|\frac{1}{N^2}\sum_{n, m=1}^N
\mathbf{1}_A(T_1^nx)\mathbf{1}_A(T_2^{n+m}x)- \E(\mathbf{1}_A,
\mathcal{I}_1)(x).\E(\mathbf{1}_A, \mathcal{I}_2)(x)\|_2 = 0.$ It
remains to show that $$\int_A \E(\mathbf{1}_A,
\mathcal{I}_1)(x).\E(\mathbf{1}_A, \mathcal{I}_2)(x)d\mu\geq
\mu(A)^3$$ if $\mathcal{I}_1\subset \mathcal{I}_2$. We have
\[\begin{aligned}
&\int_A \E(\mathbf{1}_A, \mathcal{I}_1)(x).\E(\mathbf{1}_A,
\mathcal{I}_2)(x)d\mu \\
&= \int\mathbf{1}_A(x)\E(\mathbf{1}_A,
\mathcal{I}_1)(x).\E(\mathbf{1}_A, \mathcal{I}_2)(x)d\mu =
\int\E(\mathbf{1}_A, \mathcal{I}_2)(x)\E(\mathbf{1}_A,
\mathcal{I}_1)(x).\E(\mathbf{1}_A, \mathcal{I}_2)(x)d\mu\\
&=\int \E\big[\E(\mathbf{1}_A, \mathcal{I}_2)^2,
\mathcal{I}_1\big](x)\E(\mathbf{1}_A, \mathcal{I}_1)(x)d\mu\geq
\int \E(\mathbf{1}_A, \mathcal{I}_1)^2(x)\E(\mathbf{1}_A,
\mathcal{I}_1)(x)d\mu\\
&=\int \E(\mathbf{1}_A,\mathcal{I}_1)^3d\mu\geq \big(\int
\E(\mathbf{1}_A,\mathcal{I}_1)(x)d\mu\big)^3 = \mu(A)^3
\end{aligned}
\]
This ends the proof of the corollary 1.
\subsection{Corollary 2}
For each fixed positive integer $k$ we just need to apply theorem
2 to the functions $f_i = \mathbf{1}_A$ for $1\leq i \leq
2^{k}-1$. The pointwise convergence of the averages along the
cubes of these $2^k -1$ functions to the limit $\mu(A)^{2^k}$
indicates that for $\mu$ a.e. $x$ the set
$$ \{(n_1,n_2,...,n_k)\in\Z^k: \mathbf1_A(x).\mathbf1_A(T^{n_1}_{1}x).\mathbf1_A(
  T^{n_1+n_2}_{2}x)\cdots\mathbf1_A(T^{n_1+n_2 +...+n_k}_{k}x)>\lambda \mu(A)^{2^k}\}$$
  is syndetic.

\end{document}